\documentclass{article}

\usepackage{arxiv}

\usepackage[utf8]{inputenc} 
\usepackage[T1]{fontenc}    
\usepackage{hyperref}       
\usepackage{url}            
\usepackage{booktabs}       
\usepackage{amsfonts}       
\usepackage{nicefrac}       
\usepackage{microtype}      
\usepackage{lipsum}
\usepackage{graphicx}
\graphicspath{ {./images/} }

\usepackage{amsmath,amssymb,amsfonts}
\usepackage{algorithmic}
\usepackage{textcomp}
\usepackage{nccmath}
\usepackage{soul}
\newtheorem{thm}{Theorem}

\DeclareMathOperator*{\argmin}{arg\,min}

\newcommand{\bea}{\begin{eqnarray}}
\newcommand{\eea}{\end{eqnarray}}

\newcommand{\be}{\begin{equation}}
\newcommand{\ee}{\end{equation}}

\newcommand{\ben}{\begin{equation*}}
\newcommand{\een}{\end{equation*}}

\newcommand{\bean}{\begin{eqnarray*}}
\newcommand{\eean}{\end{eqnarray*}}

\title{Accelerating Optimal Experimental Design for Robust Synchronization of Uncertain Kuramoto Oscillator Model Using Machine Learning}

\author{
 Hyun-Myung Woo\\
  Department of Electrical and Computer Engineering\\
  Texas A\&M University\\
  College Station, TX 77843, USA\\
  \texttt{larcwind@tamu.edu}\\
  \And
 Youngjoon Hong\\
  Department of Mathematics\\
  Sungkyunkwan University\\
  Suwon 16419, Republic of Korea\\
  \texttt{hongyj@skku.edu}\\
  \And
 Bongsuk Kwon\\
  Department of Mathematical Sciences\\
  Ulsan National Institute of Science and Technology\\
  Ulsan 44919, Republic of Korea\\
  \texttt{bkwon@unist.ac.kr}\\
  \And 
  Byung-Jun Yoon\\
  Department of Electrical and Computer Engineering\\
  Texas A\&M University\\
  College Station, TX 77843, USA\\
  \texttt{bjyoon@ece.tamu.edu}\\
}

\begin{document}
\maketitle
\begin{abstract}
Recent advances in objective-based uncertainty quantification (objective-UQ) have shown that such a goal-driven approach for quantifying model uncertainty is extremely useful in real-world problems that aim at achieving specific objectives based on complex uncertain systems. Central to this objective-UQ is the concept of mean objective cost of uncertainty (MOCU), which provides effective means of quantifying the impact of uncertainty on the operational goals at hand. MOCU is especially useful for optimal experimental design (OED) as the potential efficacy of an experimental (or data acquisition) campaign can be quantified by estimating the MOCU that is expected to remain after the campaign. However, MOCU-based OED tends to be computationally expensive, which limits its practical applicability. In this paper, we propose a novel machine learning (ML) scheme that can significantly accelerate MOCU computation and expedite MOCU-based experimental design. The main idea is to use an ML model to efficiently search for the optimal robust operator under model uncertainty, a necessary step for computing MOCU. We apply the proposed ML-based OED acceleration scheme to design experiments aimed at optimally enhancing the control performance of uncertain Kuramoto oscillator models. Our results show that the proposed scheme results in up to ${154}$-fold speed improvement without any degradation of the OED performance.\footnote{This work has been submitted to the IEEE for possible publication. Copyright may be transferred without notice, after which this version may no longer be accessible.}
\end{abstract}



\section{Introduction}
\label{sec:introduction}
Many real-world engineering applications involve mathematical modeling of complex systems, where the constructed models are used for designing operators--such as controllers, filters, classifiers, estimators--that can effectively achieve engineering goals of interest. For example, one may be interested in building a network model representing the transcription regulations in micro-organisms that regulate their metabolism~\cite{Niu2021}. The resulting model may be used to infer the potential impacts of modifications in the transcription regulatory network (TRN) on the metabolism of interest, for example, predicting the metabolic flux changes that result from the deletion of one or more transcription factors. In this example, the engineering goal may be predicting the optimal genetic modification in the TRN that will lead to maximizing the production of a metabolite of interest. In fact, designing optimized strains of micro-organisms for ethanol overproduction~\cite{Shen2019} is an active area of research due to its implications in efficient bio-energy production.

A fundamental challenge in the aforementioned application as well as many other real-world engineering problems involving complex systems is the difficulty of accurate model construction. While one may have ample training data for model inference, the data size may nevertheless pale in comparison to the complexity of the system being modeled. Prior knowledge, if available, may also aid in improving model construction, but the final model is likely to still have substantial uncertainties. Consequently, a critical question is how one may reliably and optimally achieve the given engineering goals in the presence of model uncertainty. Furthermore, when one has the experimental budget for the acquisition of additional data or relevant knowledge (\textit{e.g.}, via hypothesis testing), how should the experimental campaigns be designed to maximize the expected ``return on investment''?

While these are fundamental problems in modern engineering with a long and rich history~\cite{Dougherty2017,Dougherty2019}, it has been recently shown that a novel Bayesian paradigm for objective-based uncertainty quantification (objective-UQ) based on the mean objective cost of uncertainty (MOCU)~\cite{Yoon2013tsp,Yoon2021multiMOCU} can effectively address the optimal design of operators and experiments for complex uncertain systems~\cite{Hong2021, Zhao2020, Dehghannasiri15bmc, Dehghannasiri15tcbb, Broumand2015pr}. The core idea underlying the MOCU-based optimal experimental design (OED) is that, when dealing with complex uncertain models, one should quantify the model uncertainty in an objective-based manner and design experiments that can reduce the uncertainty that impacts one's operational goals. By focusing on the uncertainty that matters to the operation to be performed, the experimental budget can be efficiently used for optimizing the operational performance.
To date, the efficacy of MOCU-based OED has been demonstrated in various systems, including experimental design for robust intervention in gene regulatory networks (GRNs)~\cite{Dehghannasiri15bmc,Dehghannasiri15tcbb} and that for robust synchronization of inter-coupled Kuramoto oscillators~\cite{Hong2021}.

One practical challenge that limits the potential applicability of the MOCU-based OED scheme is its high computational cost, as discussed in~\cite{Hong2021,Dehghannasiri15bmc}. The computation of MOCU involves identifying the optimal robust operator for an uncertainty class that consists of all possible models (\textit{e.g.}, models with different parameter values) as well as evaluating expectations based on high-dimensional prior (or posterior) probability distributions. Except for very simple cases, there is no closed-form expression for the optimal robust operator and the expectations have to be evaluated numerically~\cite{Hong2021}. As a result, the evaluation of MOCU involves costly optimization to find the optimal robust operator as well as extensive sampling of the uncertain model parameters from the uncertainty class to obtain reliable estimates, which may make the cost of MOCU computation formidably high in many applications.

In this paper, we tackle this issue by adopting a machine learning (ML) approach for an efficient design of the optimal robust operator, thereby significantly accelerating the computation of MOCU as well as the MOCU-based experimental design. To the best of our knowledge, this is the first study that investigates adopting ML to accelerate MOCU-based OED. In order to develop and validate this ML-based OED acceleration scheme, we focus on designing experiments that can enhance the robust control of uncertain Kuramoto models that was investigated recently in~\cite{Hong2021}. A Kuramoto model~\cite{kuramoto1975self} consists of a network of interconnected oscillators, whose dynamics are described by coupled ordinary differential equations (ODEs). The Kuramoto oscillator model has been widely studied in various fields across engineering, physics, chemistry, and biology, due to its capability to model interesting collective behavior (\textit{e.g.}, global/partial synchronization) that emerge in complex networks~\cite{Wiesenfeld_1998,N_da_2000,hammond2007pathological,Kitzbichler2009,Breakspear_2010,Bhowmik_2012,simpson2012droop,Fernandez_2015,Skardal2015,Mohseni2017,Choi2019,guo2021overviews}. For example, a \textit{microgrid} system with droop-controlled inverters can be mathematically cast as a Kuramoto model, where the synchronization failure of the model corresponds to a power outage in the microgrid~\cite{simpson2012droop,Skardal2015,guo2021overviews,PRL12,PG13,Auto13,PNAS13}. Another interesting example is the application of the Kuramoto model for studying brain dynamics~\cite{hammond2007pathological,Kitzbichler2009,Mohseni2017,Choi2019}, where the synchronization phenomena may be associated with neurodegenerative diseases~\cite{Mohseni2017,lehnertz2009synchronization}. We show that our proposed ML-based OED acceleration scheme can improve the speed of MOCU-based experimental design by ${104\sim154}$ times without degrading the OED performance. 

The two major contributions of this paper are as follows. First, we propose an ML-based scheme for the acceleration of MOCU-based OED, which leads to significant speed improvement without performance degradation. Second, we present a comprehensive analysis of ML-based MOCU estimation and validate its performance in the context of OED.

The paper is organized as follows. In Sec.~\ref{sec:KM}, we provide a brief review of the Kuramoto model and the MOCU-based OED strategy for uncertain Kuramoto models. We propose the ML-based OED acceleration strategy in Sec.~\ref{sec:accelerate_oed}. In Sec.~\ref{sec:mocu_ml} and Sec.~\ref{sec:oed_ml}, we evaluate the performance of the proposed scheme for approximate MOCU computation and experimental design, respectively. We conclude the paper in Sec.~\ref{sec:conclusion} with further discussions and potential future research directions.


\section{Overview of Optimal Experimental Design Strategy for the Uncertain Kuramoto Model}
\label{sec:KM}
In this section, we provide a brief review of the OED strategy for uncertain Kuramoto oscillator models, which we originally proposed in our recent work~\cite{Hong2021}. We begin the section with an introduction to the Kuramoto model, followed by a brief description of the robust synchronization problem for uncertain Kuramoto models. Given an uncertain Kuramoto model, we describe how the MOCU can be used to quantify the impact of the model uncertainty on the control synchronization performance and how the MOCU-based OED strategy can be used to effectively reduce the uncertainty that matters to the objective at hand--\textit{i.e.}, optimal robust synchronization of the Kuramoto model in the presence of uncertainty.

\subsection{Uncertain class of Kuramoto models}
Consider the Kuramoto model that consists of $N$ interacting oscillators described by the following ODEs:
\begin{equation}
    \dot{\theta}_i \left( t \right) = \omega_i + \sum^N_{j=1} a_{i, j} \sin \left(\theta_j \left( t \right) - \theta_i \left( t \right) \right),\label{kuramotoModel}
\end{equation}
for ${i = 1, 2,\dots, N}$, where ${\theta_i\left(t\right)}$ is the instantaneous phase of the ${i}$th oscillator at time $t$, ${\omega_i}$ is the natural frequency of the $i$th oscillator, and ${a_{i, j}}$ is the coupling strength between the ${i}$th and ${j}$th oscillators. Kuramoto models have been widely studied to investigate the synchronization phenomena in various biological, chemical, or engineered oscillator systems, whose primary interest is whether the oscillators in a given Kuramoto model will get frequency synchronized as follows:
\begin{equation}
\lim_{t \to \infty} | \dot{\theta}_i \left(t\right) - \dot{\theta}_j \left(t\right) | = 0,\label{criterion}
\end{equation}
for ${1 \le i, j \le N}$. For example, it has been shown that modern smart grid networks referred to as microgrids can be modeled as a network of Kuramoto model oscillators, where the synchronization phenomena of the Kuramoto model are closely tied with the stability of the power grid network~\cite{PRL12,PG13,Auto13,PNAS13}. Furthermore, in neuroscience studies, brain network synchronization has been shown to be associated with various neurological disorders, where excessive neuronal activities can be represented as a global synchronization of the Kuramoto model~\cite{hammond2007pathological,Mohseni2017,Choi2019,lehnertz2009synchronization,PGR2021}. While conditions for synchronization have been extensively studied for homogeneous Kuramoto models with uniform coupling strength~\cite{ARE2008,RODRIGUES20161,ABRS}, there is yet no closed-form solution that can be used to predict the asymptotic synchronization of a general heterogeneous Kuramoto model based on its parameters.

In a real-world setting, the parameters of the Kuramoto model, which represents a complex network of oscillators, may not be completely known. For example, while it may be relatively easy to accurately estimate the natural frequency of each oscillator, in the absence of interactions with other oscillators, it will be practically challenging to accurately measure the coupling strengths between all oscillators in a large network. This uncertainty gives rise to an uncertainty class of Kuramoto models, which contains all possible Kuramoto models that are consistent with our prior knowledge regarding the true model and/or available observation data. Under this setting, our primary interest is how we can apply robust control to the uncertain Kuramoto model, comprised of a network of oscillators whose natural frequency ${\omega_{i}}$ is known but their coupling strength ${a_{i, j}}$ is only known up to a range ${a_{i, j}} \in {\left[a_{i, j}^{L}, a_{i, j}^{U}\right]}$. We denote the uncertainty class of all possible Kuramoto models as ${\mathbf{\mathcal{A}}}$, which consists of all parameter vector ${\mathbf{a} = [a_{1, 2}, a_{1, 3}, \dots, a_{N-1, N}]^T} \in {\mathbf{\mathcal{A}}}$ that satisfies the given constraints. As in the previous study~\cite{Hong2021}, we assume a prior distribution ${P_{\mathbf{\mathcal{A}}} \left( \mathbf{a} \right)}$ is uniformly distributed. However, this is not necessary. Non-uniform priors may be assumed, or custom priors may be constructed based on available prior domain knowledge~\cite{Boluki2017bmc,Boluki2017tcbb}.

\subsection{Robust control of uncertain Kuramoto models}
\label{robustControlOfUncertainKuramotoModels}
\begin{figure*}[!t]
\centering
\includegraphics[width=1\linewidth]{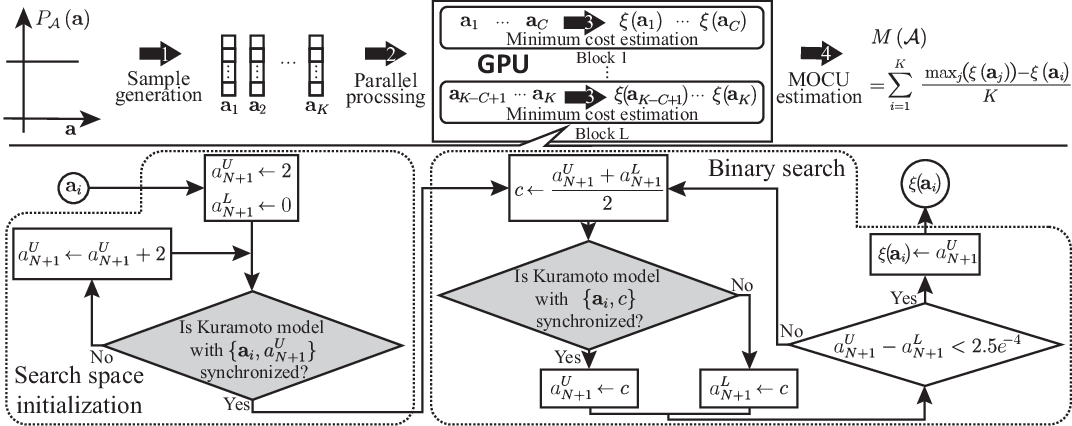}
\caption{Illustration of the original sampling-based mean objective cost of uncertainty (MOCU) computation scheme in~\cite{Hong2021}. For reliable MOCU estimation, a relatively large sample size $K$ is needed (step $1$). The sampling-based estimation scheme takes advantage of GPU programming for massive parallelization of the sampling operation. Specifically, we group the $K$ sample points ${\left\{\mathbf{a}_i\right\}}$, ${i = 1, 2, \dots, K}$, into $L$ blocks (step $2$), and the GPU processes $L$ sample points in different processing blocks in parallel (step $3$). Within block ${l}$, based on a sample point ${\mathbf{a}_i}$ that specifies a Kuramoto model (sampled from the uncertainty class), we find a valid search range ${\left[a^L_{N+1}, a^U_{N+1} \right]}$ that contains at least one valid solution that leads to global synchronization of the Kuramoto model (left bottom part). In the next phase (right bottom part), we find the solution with the smallest cost ${\xi\left(\mathbf{a}_i\right)}$ through a binary search, reducing the search range by half in every iteration. Finally, we compute the MOCU ${M\left(\mathbf{\mathcal{A}}\right)}$ of the uncertainty class ${\mathbf{\mathcal{A}}}$ based on the $K$ estimates ${\xi\left(\mathbf{a}_i\right)}$, ${i = 1, 2, \dots, K}$, (step 4).}
\label{fig1}
\end{figure*}

Suppose that we are interested in synchronizing an uncertain Kuramoto model that consists of $N$ interacting oscillators, whose interaction strengths are only known up to a range, via external control. We adopt the synchronization method proposed in~\cite{Hong2021} that introduces an additional oscillator as a global ``synchronizer'' to the original model. Let the natural frequency of this $\left(N+1\right)$th oscillator be ${\omega_{N+1}}=\frac{1}{N}\Sigma_{i=1}^N {\omega_i}$, and we assume that this control oscillator interacts with all oscillators in the original model with a uniform coupling strength $a_{i,N+1} = {a_{N+1}}$, $\forall i$, which is a control parameter. The addition of the control oscillator augments the Kuramoto model as follows:
\begin{equation}\label{kuramotoModelUnderControl}
\dot \theta_i(t) = \omega_i + \sum_{j=1}^{N} a_{i,j} \sin (\theta_j(t)- \theta_i(t)) + a_{N+1} \sin (\theta_{N+1}(t)- \theta_i(t)),
\end{equation}
for ${i=1,2,\dots,N+1}$. As the increase of the coupling strength ${a_{N+1}}$ will in practice lead to an increase of the control cost, our control objective is to find a minimum ${a_{N+1}}$ that guarantees the asymptotic frequency synchronization of the Kuramoto model despite the uncertainty. If we had complete knowledge about the coupling strength ${\mathbf{a}}$, we would be able to find the optimal (minimum) coupling strength $a_{N+1} = {\xi\left(\mathbf{a}\right)}$ that ensures synchronization by gradually increasing the value of $a_{N+1}$ from $0$ until synchronization is achieved. A more efficient approach will be to perform a binary search as illustrated in Fig.~\ref{fig1} (see the blow-up figure at the bottom). In the presence of uncertainty, we have to ensure that the control oscillator will be able to achieve synchronization for any $\mathbf{a} \in {\mathbf{\mathcal{A}}}$. For this reason, we have chosen $a_{N+1} = {\xi^\ast\left(\mathbf{\mathcal{A}}\right)}$ as follows:
\begin{equation}
\xi^\ast\left(\mathbf{\mathcal{A}}\right)=\max_{\mathbf{a}\in\mathbf{\mathcal{A}}} \xi\left(\mathbf{a}\right),    \label{eq:robust}
\end{equation}
which is the smallest $a_{N+1}$ that guarantees global synchronization of the uncertain Kuramoto oscillators.
\subsection{Objective-based quantification of the uncertainty in the Kuramoto model}
Given an uncertain Kuramoto model, the expected impact of this model uncertainty on the operational goal--in this case, the global frequency synchronization of the Kuramoto oscillators--can be quantified by the MOCU~\cite{Yoon2013tsp}. For a given uncertainty class ${\mathbf{\mathcal{A}}}$, MOCU ${M\left(\mathbf{\mathcal{A}}\right)}$ can be computed by:
\begin{equation}
    M\left(\mathbf{\mathcal{A}}\right) = E_{\mathbf{\mathcal{A}}}\left[ \xi^{\ast} \left(\mathbf{\mathcal{A}}\right) - \xi \left(\mathbf{a}\right) \right],
    \label{MOCU}
\end{equation}
where ${\xi^{\ast} \left(\mathbf{\mathcal{A}}\right)}$ is the cost of the optimal robust control and ${\xi\left(\mathbf{a}\right)}$ is the cost of the optimal control for the specific model configured with a specific parameter set ${\mathbf{a}}$. As shown in~\eqref{MOCU}, MOCU ${M\left(\mathbf{\mathcal{A}}\right)}$ quantifies the expected cost increase for applying the optimal robust control (which is inevitable to maintain robust control performance in the presence of uncertainty) instead of the model-specific optimal control (which cannot be applied in practice as the true model is unknown). In this study, the optimal robust interaction strength (cost) $\xi^{\ast} \left(\mathbf{\mathcal{A}}\right)$, ensuring that the uncertain Kuramoto model is synchronized by the added control oscillator while keeping the control cost minimum, is given by~\eqref{eq:robust}. 

\subsection{Numerical computation of mean objective cost of uncertainty}
\label{samplingBasedMOCU}
In general, there is no closed-form expression of~\eqref{MOCU}, as a result of which the MOCU ${M\left(\mathbf{\mathcal{A}}\right)}$ for uncertainty class ${\mathbf{\mathcal{A}}}$ computation requires a numerical approximation. One practical way to compute the MOCU ${M\left(\mathbf{\mathcal{A}}\right)}$ is to take a sampling-based approach to approximate it through the empirical expectation of the differential cost based on samples drawn from the distribution ${P_{\mathbf{\mathcal{A}}}\left( \mathbf{a} \right)}$.

Figure~\ref{fig1} illustrates the sampling-based MOCU computation process. First, we draw ${K}$ sample points ${\{\mathbf{a}_i\}}$, ${i= 1, 2, ..., K}$, from ${P_{\mathbf{\mathcal{A}}}\left( \mathbf{a} \right)}$. Then, for each sample point ${\mathbf{a}_i}$, which is a potential true model parameter in the uncertainty class $\mathbf{\mathcal{A}}$, we estimate the minimum coupling strength ${\xi\left( \mathbf{a}_i\right)}$ of the control oscillator that assures the asymptotic frequency synchronization of the Kuramoto model under control. To this aim, we consider a binary search to find the minimum coupling strength ${\xi\left( \mathbf{a}_i\right)}$ efficiently, as depicted in the dotted box at the bottom of Fig.~\ref{fig1}. Specifically, we start with a broad search space that contains at least one coupling strength synchronizing the system. At each iteration, we solve the ODEs of the Kuramoto model augmented with the control oscillator whose coupling strength ${a_{N+1}}$ is set to the median value ${c}$ of the current search space: ${a_{N+1}\leftarrow c = \left({a^U_{N+1}}+{a^L_{N+1}}\right)/2}$. If the system under control is synchronized, we update the upper bound of the search space to the median value: ${a^U_{N+1}\leftarrow c}$. Otherwise, we set the lower bound of the search space to the median value: ${a^L_{N+1}\leftarrow c}$. The binary search continues until we find the minimum coupling strength ${\xi\left( \mathbf{a}_i\right)}$, for the given sample point ${\mathbf{a}_i}$, which is within a specified tolerance level (set to $2.5 \times 10^{-4}$ in this study). Based on the $K$ sample points, we can obtain the MOCU ${M\left(\mathbf{\mathcal{A}}\right)}$ as follows:
\begin{equation}
    M\left( \mathbf{\mathcal{A}}\right) =  \frac{1}{K} \sum^K_{i = 1} \left( \max_j  \left(\xi\left(\mathbf{a}_j \right) \right) - \xi\left(\mathbf{a}_i \right) \right).
\end{equation}
Note that the accuracy of this numerical approximation of MOCU is dependent on the sample size ${K}$. In general, a larger $K$ generally leads to a more accurate MOCU estimation. However, at the same time, the computational cost increases as the sample size increases. We can reduce the computational time for numerical MOCU computation by exploiting parallelism. For example, estimating the optimal cost ${\xi\left( \mathbf{a}_i\right)}$ of a sample point ${\mathbf{a}_i}$ is an independent process to those of the other samples ${\mathbf{a}_j}$, ${j\neq i}$, which can be processed in a parallel manner with powerful parallel processors. In fact, the sampling-based MOCU computation in~\cite{Hong2021} takes advantage of GPU programming with Compute uniﬁed device architecture (CUDA), in which ${200}$ sample points are processed in parallel at a given time--\textit{i.e.}, ${L=200}$. However, for each sample point ${\mathbf{a}_i}$, the estimation of the minimum cost ${\xi\left( \mathbf{a}_i\right)}$ via binary search (step 3 in Fig.~\ref{fig1}) is a highly sequential process--which involves repeatedly solving the ODEs of the corresponding Kuramoto model and verifying whether or not the model is globally synchronized (\textit{i.e.}, not amenable to parallelization). In Sec.~\ref{sec:accelerate_oed}, we present a novel solution via ML that can effectively address this performance bottleneck, and thereby accelerate the numerical computational of MOCU by several orders of magnitude.

\subsection{Designing optimal experiments for effective uncertainty reduction}

The significance of objective-UQ using MOCU is that it enables the design of experiments that focus on reducing the model uncertainty that matters. More specifically, as MOCU quantifies the expected cost increase (relevant to our operational goal) due to model uncertainty, it can be used to quantify the expected impact of a potential experiment on reducing the model uncertainty that affects the operational performance, hence how effective the experiment will be in reducing the operational cost. 

The MOCU-based OED strategy for uncertain Kuramoto models has been recently proposed in~\cite{Hong2021}. In this study, a realistic experimental design space was considered, where an experiment corresponds to selecting a pair $\left(i,j\right)$ of oscillators and observing whether they get spontaneously synchronized in isolation of other oscillators and in the absence of external control. The experimental outcome is binary--either synchronized or non-synchronized--based on which the uncertainty of the coupling strength $a_{i,j} \in [a^L_{i, j}, a^U_{i, j}]$ can be reduced. Theorem~\ref{theo:sync} in~\cite{Hong2021} reproduced below gives us the necessary and sufficient condition for an oscillator pair to be frequency synchronized:
\begin{thm} \label{theo:sync}
Consider the Kuramoto model of two-oscillators:
\begin{equation}\label{Ku2}
\begin{split}
\dot \theta_1 \left(t\right)&= \omega_1 + 0.5a \sin\left(\theta_2\left(t\right)- \theta_1\left(t\right)\right),\\
\dot \theta_2 \left(t\right)&= \omega_2 + 0.5a \sin\left(\theta_1\left(t\right) - \theta_2\left(t\right)\right),
\end{split}
\end{equation}
with the initial angles $\theta_1\left(0\right), \theta_2\left(0\right)\in[0,2\pi)$.
Then, for any solutions  $\theta_1 \left(t\right)$ and $\theta_2\left(t\right)$ to~\eqref{Ku2}, there holds 
$|\dot\theta_1(t) - \dot\theta_2(t) | \to 0$ as $t\to\infty$ 
if and only if 
$|\omega_1-\omega_2|\le a$. \hspace{3.6in} $\blacksquare$
\end{thm}
According to Theorem~\ref{theo:sync}, the Kuramoto oscillator pair $\left(i,j\right)$ becomes frequency synchronized $\lim_{t \to \infty}\lvert \dot{\theta}_i \left(t\right) - \dot{\theta}_j \left(t\right)\rvert = 0$ if and only if $\frac{\left\lvert\omega_i - \omega_j \right\rvert}{2} \leq a_{i,j}$. As a result, if the two oscillators are observed to be synchronized, we can decrease the upper bound ${a^U_{i, j}}$ to ${{\left\lvert\omega_i - \omega_j \right\rvert}/{2}}$. Otherwise, we can increase the lower bound ${a^L_{i, j}}$ to ${{\left\lvert\omega_i - \omega_j \right\rvert}/{2}}$. Since the experimental outcome is unknown in advance, we need to consider both possible outcomes to quantify the expected impact of a given experiment on reducing the objective uncertainty. To formalize this, let ${O_{i, j}}$ be a binary random variable representing the outcome of the pairwise synchronization experiment for the oscillator pair $\left(i,j\right)$. Then, the expected remaining MOCU ${R\left(i, j\right)}$ is given by:
\begin{equation}
\begin{split}
R\left(i, j\right) &= E_{O_{i, j}}[M\left(\mathbf{\mathcal{A}} | O_{i, j}\right)]\\
&=\sum_{o \in \left\{ 0, 1\right\}} P\left(O_{i, j} = o\right) M\left(\mathbf{\mathcal{A}} | O_{i, j} = o \right),\label{erMOCU}
\end{split}
\end{equation}
where $M\left(\mathbf{\mathcal{A}} | O_{i, j}\right)$ is the conditional MOCU given $O_{i, j}$. The conditional MOCU $M\left(\mathbf{\mathcal{A}} | O_{i, j} = o \right)$ given an experimental outcome $O_{i, j} = o$ can be computed by reducing the uncertainty class as previously described and numerically computing the MOCU of this reduced uncertainty class. The probability $P\left(O_{i, j} = o\right)$ can be derived in a straightforward manner, based on ${P_{\mathbf{\mathcal{A}}}\left( \mathbf{a} \right)}$ (see \cite{Hong2021} for further details). The ${R\left(i, j\right)}$ in~\eqref{erMOCU} quantifies the MOCU that is expected to remain after performing the pairwise synchronization experiment for the pair $\left(i,j\right)$.

So, how should we prioritize the potential $N \choose 2$ experiments? Naturally, the optimal choice will be to choose the experiment with the smallest ${R\left(i, j\right)}$:
\be
    \left(i^\ast, j^\ast\right)=\argmin_{\left(i,j\right)\in\mathbf{\mathcal{E}}} R\left(i,j\right), \label{eq:oed}
\ee
as experiment $\left(i^*,j^*\right)$ is expected to most effectively reduce the objective uncertainty among all potential experiments. In practice, rather than performing a single best experiment, we may perform a sequence of experiments prioritized by~\eqref{eq:oed}. In theory, $R\left(i,j\right)$ needs to be re-estimated after performing the predicted optimal experiment and observing its outcome, as it changes the uncertainty class, hence the expected remaining MOCU for the potential subsequent experiments. However, empirically, $R\left(i,j\right)$ computed based on the original uncertainty class $\mathbf{\mathcal{A}}$ is a robust indicator of the efficacy of the potential experiments, which we will demonstrate in Sec.~\ref{sec:oed_ml}.

\subsection{Computational complexity of experimental design}
The overall computational complexity for predicting the optimal experiment is as follows: 
\be
    O(TKN^4L^{-1} \log \epsilon ), \label{eq:oed_complexity}
\ee
where $T$ is the time duration for solving the ODEs using the Runge-Kutta method (to check for asymptotic global frequency synchronization among the Kuramoto model oscillators), $K$ is the sample size for numerical computation of MOCU, $N$ is the number of oscillators in the Kuramoto model, ${L}$ is the number of parallel processing blocks in GPU, and $\epsilon$ is the tolerance level for the binary search (set to $\epsilon = 2.5 \times 10^{-4}$ in this study). Note that the complexity for computing MOCU is $O(TKN^2L^{-1} \log \epsilon)$, where predicting the optimal experiment involves computing MOCU $2 \cdot {N \choose 2}$ times to calculate $R\left(i,j\right)$ given by~\eqref{erMOCU} for all oscillators pairs. As we can see in~\eqref{eq:oed_complexity}, the computational cost for OED sharply increases as the size $N$ of the Kuramoto model increases, which limits the practical applicability of the OED scheme~\cite{Hong2021} for large models. For example, when ${T=5}$, ${K=20,480}$, and ${L=128}$, respectively, identifying the optimal experiment $\left(i^*,j^*\right)$ for the uncertain Kuramoto model operating on five oscillators required ${650}$ seconds on average. However, it took ${3,171}$ seconds to determine the optimal experiment $\left(i^*,j^*\right)$ for the uncertain Kuramoto model with seven oscillators.

\section{Accelerating Experimental Design via Machine Learning}
\label{sec:accelerate_oed}
We propose an ML approach for accelerating the quantification of the objective system uncertainty. As we discussed in the previous section, in real-world applications that typically involve the control consisting of highly non-linear sequential operations, the effective computational complexity is critically dependent on the computational complexity of the control rather than the number of samples. The proposed approach learns a surrogate model for (part of) the operations of the control for estimation of the control cost for a system, thereby reducing the effective computational complexity that cannot be further reduced by parallelism. Recently, there has been an increasing number of studies investigating the application of deep learning (DL) methods to scientific computation, including approximating and solving differential equations (DEs) (\textit{e.g.}, see~\cite{L712178, Han8505, RAISSI2019686} and references therein). However, it is worth noting that the primary focus of our current study does not lie in solving ODE systems via deep network models but in the accelerated design of optimal experiments based on the objective-based UQ via the concept of MOCU. Rather than aiming at a fast solution of DEs, our goal is to efficiently design experiments that can most effectively reduce model uncertainty, thereby optimally enhancing the control performance of uncertain Kuramoto oscillator models.

As discussed in Sec.~\ref{samplingBasedMOCU}, estimating the MOCU of the uncertain Kuramoto model based on the sampling approach involves a binary search for each sample ${\mathbf{a}_i}$, where at each iteration solving the corresponding ODEs and determining if the system under control is synchronized or not. From a broad perspective, at each iteration, these operations, the gray box in Fig.~\ref{fig1}, are nothing but a binary classification problem. Hence, if we collect enough samples to build an accurate classifier, we replace such a process with the binary classifier, which is computationally efficient. In this study, we considered a fully connected neural network (fcNN) with only one hidden layer, possibly the simplest ML structure that we can think of.

\begin{figure*}
\centering
\includegraphics[width=1\linewidth]{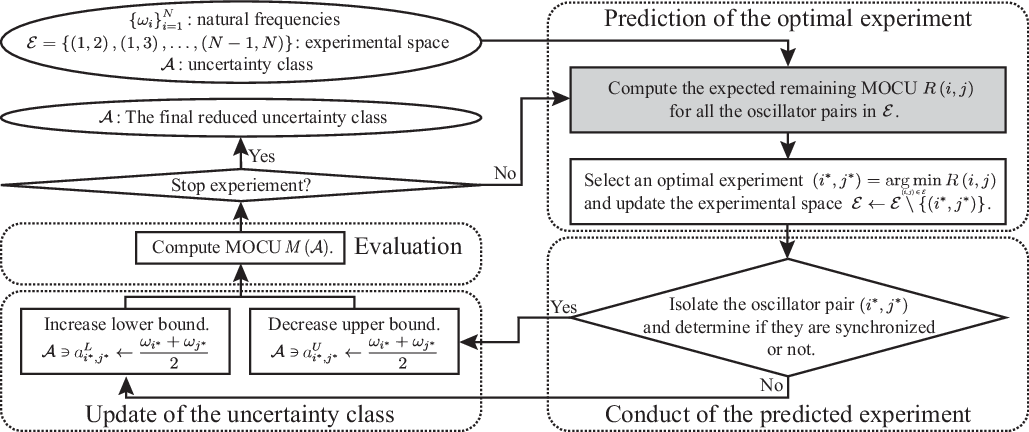}
\caption{Illustration of the MOCU-based optimal experimental design (OED) loop. First, we compute the expected remaining MOCU ${R \left( i,j\right)}$ for all possible experiments ${\left(i, j\right)}$ in the experimental design space $\mathbf{\mathcal{E}}$ based on the current uncertainty class $\mathbf{\mathcal{A}}$. Next, we identify the optimal experiment ${\left(i^\ast, j^\ast\right)}$ that has the smallest expected remaining MOCU such that ${\left(i^\ast, j^\ast\right)=\argmin_{\left(i,j\right)\in\mathbf{\mathcal{E}}} R\left(i,j\right)}$. In the second phase (right bottom), we conduct the selected experiment ${\left(i^\ast, j^\ast\right)}$ and remove the performed experiment from the experimental space $\mathbf{\mathcal{E}}$. Specifically, in this experiment, we isolate the selected oscillator pair ${\left(i^\ast, j^\ast\right)}$ and determine whether or not they get synchronized without external control. Based on the experimental outcome, we update the uncertainty class accordingly~\cite{Hong2021}. Finally, we evaluate the actual efficacy of the conducted experiment by computing the MOCU of the updated uncertainty class $\mathbf{\mathcal{A}}$. We iterate this experimental loop until the experimental space becomes empty (\textit{i.e.}, there are no more experiments left to be performed).}
\label{fig2}
\end{figure*}

\begin{figure*}
\centering
\includegraphics[width=1\linewidth]{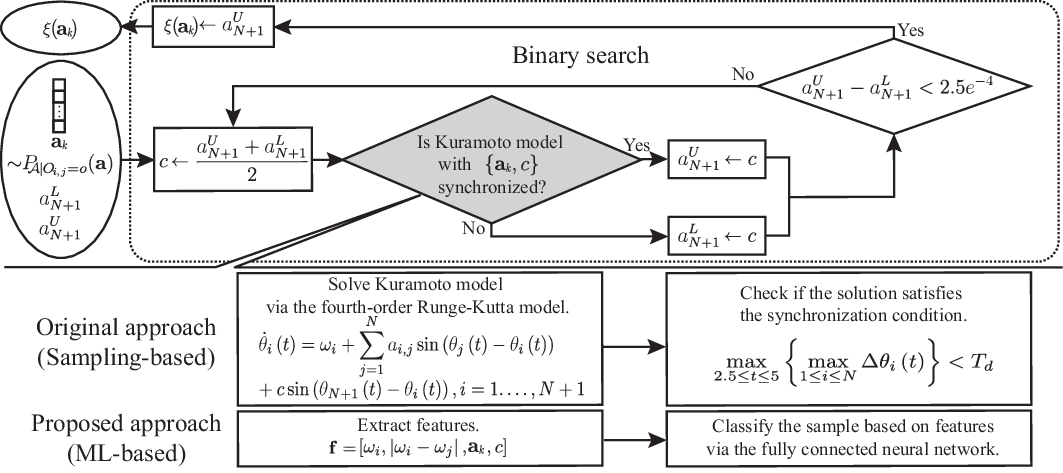}
\caption{Comparison between the original sampling-based estimation scheme adopted in~\cite{Hong2021} and the proposed machine learning-based (ML-based) estimation scheme. The proposed scheme obviates the need for repeatedly solving the coupled ordinary differential equations (ODEs) within the binary search routine to find the optimal robust coupling strength illustrated in Fig.~\ref{fig1}. This significantly enhances the computational efficiency of MOCU estimation.}
\label{fig3}
\end{figure*}

The proposed approach on the MOCU-based OED framework is realized by replacing part of the operations of control with the trained model for the estimation of the expected remaining MOCU ${R\left(i, j\right)}$ highlighted in gray in Fig.~\ref{fig2}. Hence we focus on the difference in quantifying the expected remaining MOCU ${R\left(i, j\right)}$ between the proposed ML-based approach and the original approach that manually determines the synchronization of the Kuramoto model. To compute the expected remaining MOCU ${R\left(i, j\right)}$ we first need to estimate conditional MOCU ${M\left(\mathbf{\mathcal{A}} | O_{i, j} = o \right)}$ given the experimental outcome ${O_{i, j} = o \in \{0, 1\}}$ (\textit{i.e.}, synchronized or not) as derived in~\eqref{erMOCU}. Specifically, as we described in Sec.~\ref{samplingBasedMOCU}, we compute the control cost ${\xi\left(\mathbf{a}_i\right)}$ of all samples ${\mathbf{a}_i}$, ${i = 1, 2, \dots, K}$, drawn from the posterior uncertainty class distribution ${P_{\mathbf{\mathcal{A}}| O_{i, j} = o} \left( \mathbf{a} \right)}$ updated according to the experimental outcome ${O_{i, j} = o}$ as shown in Fig.~\ref{fig2}. Figure~\ref{fig3} shows the difference in estimating the control cost ${\xi\left(\mathbf{a}_i\right)}$ of sample ${\mathbf{a}_k}$ between the proposed approach and the original approach. Both approaches find a numerical solution through the binary search that is a sequential process. At each iteration, the coupling strength ${a_{N+1}}$ of the control oscillator is set to the midpoint ${c\leftarrow\left({a^U_{N+1}}+{a^L_{N+1}}\right)/2}$ of the search space. The original approach solves the Kuramoto model determined by the sample ${\mathbf{a}_k}$ and midpoint ${c}$ and determines if the solutions are synchronized or not according to criterion~\eqref{criterion}. On the other hand, the proposed ML-based approach extracts features based on the natural frequencies ${\omega_i}$, sample ${\mathbf{a}_k}$, and midpoint ${c}$ and classifies the feature vector. The search space is then halved according to the outcome. Note that the computational complexity of the original approach is critically dependent on the time precision and simulation time. Less time precision and shorter simulation time can reduce the overall computational complexity, but such parameters significantly affect the estimation accuracy of the MOCU. On the other hand, MOCU-based OED with the proposed approach is free from such a trade-off at the inference phase as the features are independent of the parameters. 

\section{Evaluation of Machine Learning-based Estimation of Approximate Mean Objective Cost of Uncertainty}
\label{sec:mocu_ml}
\subsection{Configuration of the simulation}
In this section, we demonstrate the efficacy of the proposed ML approach in accelerating the speed of objective-UQ, resulting in a very efficient OED. As described in the previous sections, we consider the OED for the Kuramoto model under uncertainty, where one's operational objective is to ensure synchronization of the model by adding an oscillator for control. For validation, we considered two experimental setups based on uncertain Kuramoto models with five and seven oscillators, respectively. As a reference ODE solver, we used the fourth-order Runge-Kutta method to solve the Kuramoto model sampled at the sampling frequency ${f_s}$ of ${160}$Hz for five seconds. To determine whether the Kuramoto model is synchronized or not, we used the following criterion: 
\begin{equation}
    \max_{2.5\leq t\leq 5} ( \max_{1\leq i \leq N} \Delta \theta_i \left(t\right)) < T_d,\label{syncCriterion}
\end{equation}
where ${\Delta \theta_i \left(t\right) \triangleq \theta_i \left(t + ({1}/{f_s})\right) - \theta_i \left(t\right)}$, ${\theta_i\left(t\right)}$ is the instantaneous phase of the ${i}$th oscillator, and ${T_d}$ is a threshold of tolerance. We set ${T_d}$ to ${0.001}$. To estimate the MOCU of a given uncertainty class, we randomly drew $20,480$ sample points from the uncertainty class (\textit{i.e.}, ${K = 20,480}$). We used a \textit{Lambda workstation} equipped with \textit{Intel i$9$-$9960$X}, $128$GB memory, and \textit{GeForce RTX $2080$ Ti} for the simulations.

At the core of the proposed method lies a binary classifier that accurately classifies the global frequency synchronization of the model when a control oscillator is introduced. To train an accurate classifier, we used an fcNN model with one hidden layer. In that regard, it is essential to extract representative features from the parameters that define the Kuramoto model, such as the number of oscillators, natural frequencies, initial phases, or coupling strength values between oscillators. Inspired by Theorem~\ref{theo:sync}, which gives us the necessary and sufficient condition for pairwise frequency synchronization of Kuramoto oscillators, we used the natural frequencies, the absolute difference between the frequencies, and the corresponding coupling strength values as features. More specifically, given a parameter set that fully determines the Kuramoto model operating on $N + 1$ oscillators, we first sort all the natural frequencies in descending order and rearrange the coupling strength accordingly. Then, we construct the corresponding feature set that consists of the sorted natural frequencies, the absolute difference of the natural frequencies of all oscillator pairs, and their coupling strengths. Note that this arrangement makes the feature set highly structured but does not affect the characteristics of the Kuramoto model. To accurately label a given sample point (the feature set of a given Kuramoto model), we used the fourth-order Runge-Kutta method with a much longer simulation time ${T}$ of ${400}$ seconds to determine whether the model reaches global frequency synchronization or not. Besides, we rigorously determined the synchronization of the Kuramoto model based on more stringent criteria. For the labeling purpose, we consider that a Kuramoto model is synchronized if both of the following two conditions are satisfied: First, frequencies of all oscillators rounded to the sixth decimal place are equal for the last ${20}$ $\left(T*0.95\right)$ seconds. Second, the sum of absolute change in the coherence value ${r\left(t\right)}$ of the order parameter ${r\left(t\right) e^{j \psi \left(t\right)} = \frac{1}{N} \sum^{N}_{i = 1} e^{j\theta_i\left(t\right)}}$ is less than $10^{-6}$ for the last ${20}$ seconds. Note that if the results for the two conditions differ, we excluded the sample point from the training dataset. Since the detailed structure and the training process of the classifier varies depending on the simulation scenario, we provide further details in the corresponding subsections. The source code used for the simulations whose results are presented in this study can be found at \href {https://github.com/bjyoontamu/Kuramoto-Model-OED-acceleration}{https://github.com/bjyoontamu/Kuramoto-Model-OED-acceleration}.

\subsection{Performance evaluation}
To assess the performance of the proposed approach and compare it to other existing approaches, we performed a wide range of evaluation experiments. First, we evaluated the efficacy of the shallow fcNN model, which is adopted in this study, in predicting the global synchronization of the Kuramoto model oscillators. For this purpose, we examined the asymptotic behavior of the fcNN model by estimating its prediction accuracy as a function of increasing training data size. Next, we computed the Pearson's correlation coefficient between the MOCU values computed by the proposed scheme and the original sampling-based scheme~\cite{Hong2021}, respectively. Furthermore, the time complexity of each scheme was assessed to compare efficiency. Finally, we evaluated the efficacy of the proposed scheme in predicting experiments that can effectively reduce model uncertainty. To this aim, we compared the changes in the ``objective uncertainty'' (estimated by MOCU) after conducting the experiment selected by different experimental design schemes. We also examined how accurately each scheme is able to predict the true optimal experiment, and how it affects the overall experimental performance.

\subsection{Comparison between machine learning-based and sampling-based mean objective cost of uncertainty estimation}

In order to validate the efficacy of the proposed method that incorporates ML-based predictions into MOCU estimation, we first directly compare the MOCU values from the ML-based approach and the sampling-based approach that we considered in the previous work~\cite{Hong2021}.

\subsubsection{Mean objective cost of uncertainty estimation for uncertain Kuramoto model with five oscillators}
\label{MOCU5}
As a first experimental scenario, we considered an uncertain Kuramoto model that consists of five oscillators that do not get spontaneously synchronized in the absence of external control. In this experiment, we adopted the identical experimental setup in the previous work~\cite{Hong2021} for direct comparison. Specifically, we assumed that the five oscillators have the natural frequencies of ${-2.50}$, ${-0.6667}$, ${1.1667}$, ${2.0}$, and ${5.8333}$, respectively. The natural frequency of the additional (\textit{i.e.}, $6$th) control oscillator was set to the average frequency of the five oscillators (${\omega_{6} = 1.1667}$). Besides, we set the initial phase of all the oscillators to zero. Finally, we used the uncertainty class defined as follows:

\begin{equation}\label{U5}
\mathbf{a}^U =
\left[\begin{matrix}
  1.0541 & 0.6325 & 0.7762 & 1.4375 & 1.0542 & 0.6900 & 1.6819 & 0.4791 & 2.6833 & 2.2041
\end{matrix}\right]^T,
\end{equation}
\begin{equation}\label{L5}
\mathbf{a}^L =
\left[\begin{matrix}
  0.7791 & 0.4675 & 0.5737 & 1.0625 & 0.7792 &  0.5100 & 1.2431 & 0.3541 & 1.9833 & 1.6291 
\end{matrix}\right]^T.
\end{equation}

To train the classifier, we generated $40,000$ sample points (a set of $20,000$ parameter values that result in synchronization and another set of $20,000$ parameter values that do not) from a multivariate uniform distribution whose support completely covers the range of the parameters in the uncertainty class at hand. Specifically, a parameter set has six real-values from the uniform distribution with a range of ${\left( -2\pi, 2\pi\right)}$ as natural frequencies of the six oscillators ${\omega_i}$, ${i=1,2,\dots,6}$, and ten coupling strength values ${a_{i,j}}$, ${1 \le i < j \le 6}$, between oscillators ranging from ${0.25 \left|\omega_i - \omega_j\right|}$ to ${2.35 \left|\omega_i - \omega_j\right|}$. To build the classifier, we sorted the six natural frequencies in descending order and rearranged the coupling strength values accordingly. Then, we extracted the following features: the sorted natural frequencies, the absolute difference of the natural frequencies of all oscillator pairs, and their coupling strengths. Finally, we trained an fcNN model with a single hidden layer, whose width is three times the number of features, until the model is capable of classifying all the $40,000$ sample points in the training dataset perfectly. We validated the trained model in terms of its asymptotic classification accuracy by assessing the accuracy as a function of the training data size. This result is shown in Fig.~S1 in the supplemental material.

\begin{figure}
\centering
\includegraphics{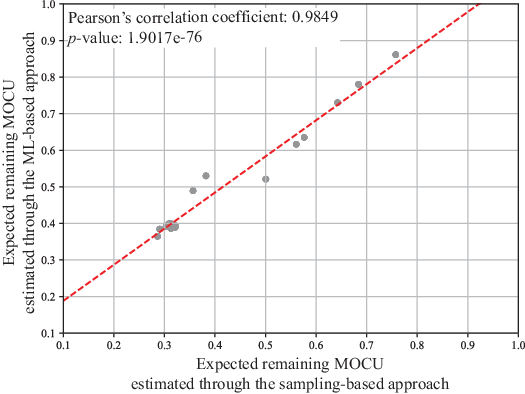}
\caption{The scatter plot shows the expected remaining MOCU values for the uncertain five-oscillator Kuramoto model estimated using the proposed ML-based approach and the original sampling-based approach in~\cite{Hong2021}. As shown, the estimated values are highly correlated to each other.}
\label{fig4}
\end{figure}

We started with the original uncertainty class defined in~\eqref{U5} and~\eqref{L5} and estimated the expected remaining MOCU of random oscillator pairs through both approaches one hundred times while randomly changing the true model (assumed to be unknown). Figure~\ref{fig4} is a scatter plot that shows the comparison between the expected remaining MOCU values computed by different methods. As shown in Fig.~\ref{fig4}, the expected remaining MOCU values computed by the proposed ML-based method and the original sampling-based method display a strong linear relationship. The Pearson's correlation coefficient was $0.9849$ with a ${p}$-value of $1.9017e^{-76}$. This plot shows that the ML-based computational scheme has the potential to effectively replace the costly sampling-based scheme without affecting the MOCU-based OED performance, as it will likely not affect the ranking of potential experiments. In terms of computational cost, the ML-based approach was able to compute the expected remaining MOCU in ${0.1110}$ seconds (on average) for a given uncertainty class, while it took ${818.7}$ seconds (on average) for the sampling-based approach. These results clearly show the advantages of the proposed approach in efficiently quantifying the objective uncertainty. 

\subsubsection{Mean objective cost of uncertainty estimation for uncertain Kuramoto model with seven oscillators}
\label{MOCU7} 
Some major disadvantages of the original sampling-based method are its high computational cost and the practical difficulty of scaling the method for large models with many oscillators. To examine the computational cost increase and scalability for larger models, we next considered an uncertainty class of Kuramoto models with seven oscillators. This increases the time for solving the ODEs and the number of possible experiments also increases from ${5\choose 2} = 10$ to ${7\choose 2} = 21$. Here we set the natural frequency of the oscillators to $-3.4600$, $-1.9611$, $-0.6754$, $-0.3806$, $-0.3675$, $6.1161$, and $8.3287$, respectively. We assumed that the natural frequency of the control oscillator (\textit{i.e.}, $8$th oscillator) is the average frequency of the seven oscillators (${\omega_{8} = 1.0857}$). We considered the uncertainty class shown below:

\begin{equation}\label{U7}
\mathbf{a}^U =
\!\begin{aligned}
&
\left[\begin{matrix}
  0.848 & 0.988 & 1.446 & 1.607 & 3.820  & 0.915 & 0.400
\end{matrix}\right.\\[-2pt]
&\ 
\left.\begin{matrix}
  0.850 & 0.419 & 4.162 & 1.090  & 0.122 & 0.039 & 2.124
\end{matrix}\right.\\[-4pt]
&\ 
\left.\begin{matrix}
  0.872 & 0.007 & 2.737 & 1.804 & 1.360  & 0.744 & 1.174
\end{matrix}\right]^T,
\end{aligned}
\end{equation}
\begin{equation}\label{L7}
\mathbf{a}^L =
\!\begin{aligned}
&
\left[\begin{matrix}
  0.073 & 0.172 & 0.153 & 0.054 & 0.501 & 0.463 & 0.043
\end{matrix}\right.\\[-2pt]
&\ 
\left.\begin{matrix}
  0.015 & 0.096 & 0.501 & 0.103 & 0.007 & 0.009 & 0.139
\end{matrix}\right.\\[-4pt]
&\ 
\left.\begin{matrix}
  0.408 & 0.000 & 0.131 & 0.119 & 0.300 & 0.286 & 0.131
\end{matrix}\right]^T.
\end{aligned}
\end{equation}
As in the previous experiment for the Kuramoto model with five oscillators, we set the initial phase of all oscillators to zero.

As the size of the parameter set is much greater for this Kuramoto model, we generated the training data in a more tailored way. Rather than generating the sample points (\textit{i.e.}, Kuramoto model parameter sets) with random natural frequencies within a specific range as we did for the five oscillator model, we fixed the natural frequencies to $-3.4600$, $-1.9611$, $-0.6754$, $-0.3806$, $-0.3675$, $6.1161$, $8.3287$, and ${1.0857}$ in this example. For the coupling strength values, we drew them from the uniform distribution for the uncertainty class, whose support is defined in~\eqref{U7} and~\eqref{L7}. In this manner, we collected $50,000$ sample points per label according to the same criteria we used for the five oscillator case. Then, we extracted the feature values as described previously for the five oscillator case and trained the classifier using an fcNN with a single hidden layer, whose width is four times the number of features. Figure~S1 in the supplemental material shows that this model quickly learns the classification boundary, where the classification accuracy rapidly converges to $100$\% as the size of the training data increases.

\begin{figure}
\centering
\includegraphics{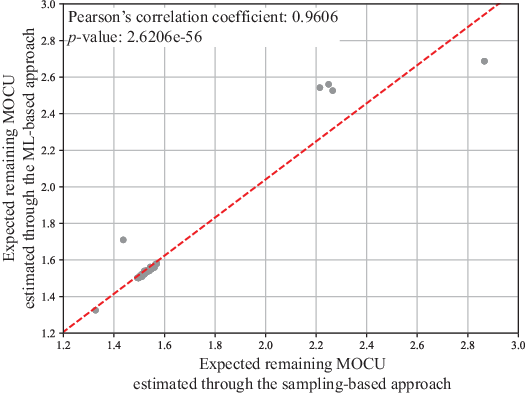}
\caption{The scatter plot shows the expected remaining MOCU values for the uncertain seven-oscillator Kuramoto model estimated using the proposed ML-based approach and the original sampling-based approach in~\cite{Hong2021}. As before, the estimated values show a high correlation.}
\label{fig5}
\end{figure}

As before, we started with the original uncertainty class defined in~\eqref{U7} and~\eqref{L7} and computed the expected remaining MOCU of random oscillator pairs using the ML-based method and the sampling-based method. We repeated this until we collected a hundred expected remaining MOCU values per method. Figure~\ref{fig5} shows the scatter plot that compares the expected remaining MOCU values computed by the two methods. Again, we can see that there is a strong linear relationship between the computed values. The Person's correlation coefficient was ${0.9606}$ with a ${p}$-value of ${2.6206e^{-56}}$. In terms of computational cost, it took ${0.6953}$ seconds (on average) for the ML-based method to compute the expected remaining MOCU, which was still less than a second although the experimental design space has grown from ${\frac{\left(5\times4\right)}{2} = 10}$ experiments to ${\frac{\left(7\times6\right)}{2} = 21}$. It took the sampling-based approach  ${3,684.9}$ seconds (on average) to compute the expected remaining MOCU values, which shows that our proposed method makes the computation $5,298$ times faster at practically identical accuracy. These results clearly show the advantages of the proposed ML-based approach in quantifying the objective model uncertainty.

\section{Performance of Experimental Design Using Machine Learning-based Mean Objective Cost of Uncertainty Estimation}
\label{sec:oed_ml}
We compared the OED performance of the proposed ML-based method against three existing approaches:
\begin{itemize}
    \item \textbf{Sampling-based approach}: the original approach proposed in~\cite{Hong2021} based on the MOCU framework, where a fourth-order Runge-Kutta method is to solve the Kuramoto model to determine synchronization.
    \item \textbf{Entropy-based approach}: the experiment is chosen for the oscillator pair whose coupling strength value has the largest entropy to reduce this uncertainty.
    \item \textbf{Random approach}: the experiment is randomly selected from the experimental design space. 
\end{itemize}
For the MOCU-based OED schemes (\textit{i.e.}, ML-based and sampling-based computations), we consider the following OED strategies. In the first approach (marked as \textit{iterative} in the figures), we re-estimate the expected remaining MOCU for the remaining experiments in each iteration, after performing the predicted optimal experiment and updating the uncertainty class based on the observed experimental outcome. In the second approach, we estimate the expected remaining MOCU only based on the initial uncertainty class and prioritize all experiments based on this result. While this approach is theoretically suboptimal, it significantly reduces the overall computational cost and empirically shows comparable performance to the iterative scheme, as we will show in this section.

\subsubsection{Optimal experimental design for uncertain Kuramoto models with five oscillators}
\begin{figure}
\centering
\includegraphics{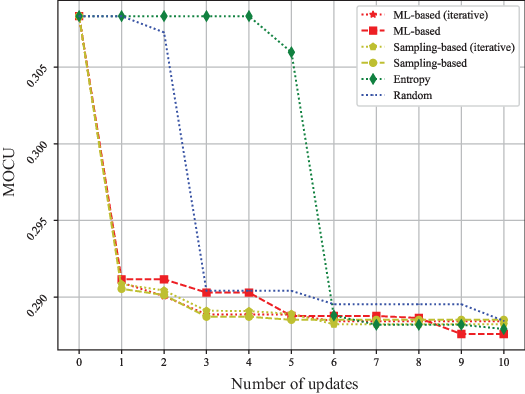}
\caption{Performance comparison of various experimental design strategies for the uncertain five-oscillator Kuramoto model considered in~\cite{Hong2021}. The results show that the three MOCU-based OED schemes perform similarly, regardless of how MOCU was estimated. The MOCU-based schemes clearly outperform other schemes as reported in~\cite{Hong2021}.}
\label{fig6}
\end{figure}

First, we conducted OED simulations for the same five-oscillator Kuramoto model considered in the previous study~\cite{Hong2021} for direct comparison. We used identical model parameters described in Sec.~\ref{MOCU5}. The true (unknown) model ${\mathbf{a}}$ was assumed to be as follows:
\begin{equation}
\mathbf{a} =
\left[\begin{matrix}
  0.9166 & 0.55 & 0.675 & 1.25 & 0.9167 & 0.6 & 1.4625 & 0.4166 & 2.3333 & 1.9166 
\end{matrix}\right]^T.
\end{equation}
Figure~\ref{fig6} shows the experimental design performance of the different algorithms, where the objective uncertainty (quantified by MOCU) is shown as a function of the number of experimental updates (iterations). As shown in Fig.~\ref{fig6}, the proposed ML-based approach with iterative re-estimation (red dotted line with asterisks) showed the nearly identical performance to sampling-based methods (both iterative and non-iterative schemes, shown in yellow lines). All three schemes reached the near minimum MOCU within only three experimental updates. The non-iterative ML-based scheme (red dashed line with squares) also identified the first optimal experiment accurately and showed comparable performance in the later updates with the other three MOCU-based OED schemes. All four MOCU-based OED schemes (both ML-based and sampling-based) significantly outperformed the entropy-based and random approaches, resulting in much sharper uncertainty reduction within fewer experimental updates.

\begin{figure}
\centering
\includegraphics{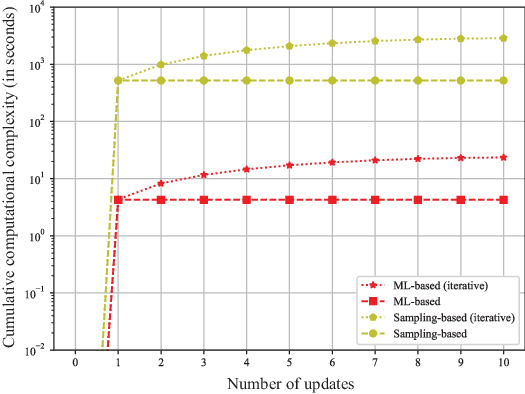}
\caption{Cumulative computational cost for experimental design. As shown, the proposed ML-based estimation clearly outperforms the original sampling-based estimation~\cite{Hong2021} in terms of efficiency, where their costs differ by two orders of magnitude. For both ML-based/sampling-based schemes, iterative estimation requires further computations, as the uncertainty class is updated after each experiment, based on which the remaining expected MOCU values are assessed again.}
\label{fig7}
\end{figure}

Figure~\ref{fig7} compares the overall computational cost between the ML-based OED schemes and the sampling-based OED schemes. The entropy-based approach and the random approach are not shown, as their computational cost is fixed and negligible. As we can see in Fig.~\ref{fig7}, the proposed ML-based OED approaches, marked as red, showed significantly lower time complexity compared to the sampling-based OED approaches. Note that the ML-based methods (red dotted lines) were significantly faster compared to the sampling-based methods, despite maintaining equivalent OED performance.

\begin{figure}
\centering
\includegraphics{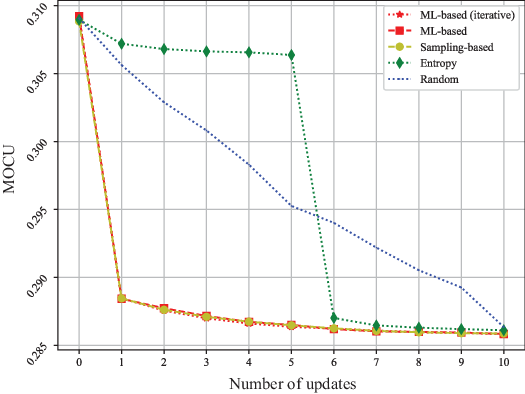}
\caption{Average performance of various experimental design strategies for uncertain Kuramoto models with five oscillators. The experiments have been repeated one hundred times by sampling potential true models from the uncertainty class. As shown, all three MOCU-based methods lead to the best performance. Random selection results in linear uncertainty reduction as expected.}
\label{fig8}
\end{figure}

\begin{figure}
\centering
\includegraphics{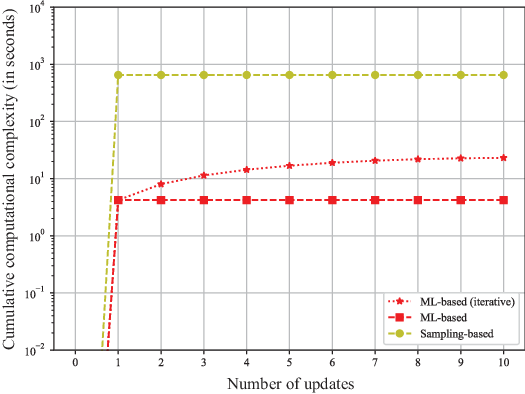}
\caption{Average cumulative computational cost (in seconds) for the different experimental design strategies for uncertain Kuramoto models with five oscillators.}
\label{fig9}
\end{figure}

Next, we repeated the experiment based on one hundred different true models randomly drawn from the uncertainty class (\textit{i.e.}, different coupling strength values were drawn from the prior distribution of the uncertainty class). The results of these large-scale experiments are shown in Fig.~\ref{fig8} and Fig.~\ref{fig9}. Note that we excluded the iterative sampling-based OED method due to its excessive requirement of computational time. As shown in these figures, the proposed ML-based method without iterative re-estimation of the expected remaining MOCU showed identical performance to other best performers. Random experimental selection (blue dotted line) yielded a linearly decreasing MOCU curve, as we would expect on average. The entropy-based method showed similar performance as before (see Fig.~\ref{fig6}). Computational cost in Fig.~\ref{fig9} shows a similar trend as before (see Fig.~\ref{fig7}). Furthermore, Fig.~S2 in the supplemental material shows the RainCloud plot~\cite{allen2019raincloud} that depicts the instantaneous performance of the different methods measured in terms of the remaining uncertainty (measured by MOCU) after performing the first experiment selected by the respective methods. As we can see from Fig.~S2, all three MOCU-based OED schemes consistently yield the best overall performance.

\begin{figure}
\centering
\includegraphics{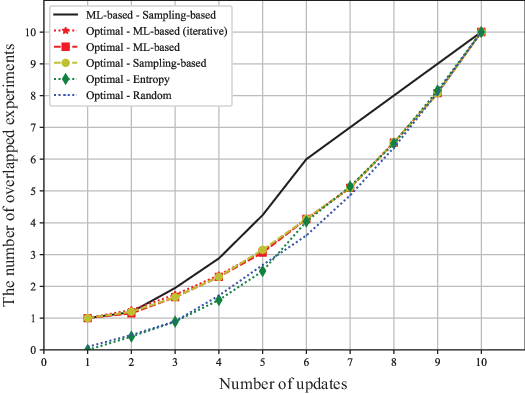}
\caption{Comparison between the optimal sequence of experiments predicted by different OED strategies for uncertain Kuramoto models with seven oscillators. The $y$-axis shows the number of common experiments within the first $k$ experiments predicted by two different methods.}
\label{fig10}
\end{figure}

Finally, we compared the experimental sequences determined by the ML-based method and the sampling-based methods, respectively, to further investigate if the proposed ML approach can practically replace the sampling-based method for prioritizing the experiments in the experimental design space. The vertical axis corresponds to the number of intersecting experiments in the first $k$ experiments predicted by two different methods. If two methods predict the identical experimental sequence, the resulting curve will be a straight line (with unit slope). For example, the black line in Fig.~\ref{fig10} compares the ML-based method and the sampling-based method. From Fig.~\ref{fig10} we can see that the proposed ML-based method (without re-estimation) always identified the same first experiment as the sampling-based method in all one hundred evaluations. By comparing the true optimal experimental sequence (\textit{i.e.}, predicted by an ``oracle'') and the sequences predicted by the ML-based method, we can see that the first optimal experiment was always accurately predicted. In fact, results in Fig.~\ref{fig8} show that the first experiment leads to the most significant drop in model uncertainty, and all MOCU-based OED schemes (both ML-based and sampling-based) accurately predict this critical experiment. We also note that the entropy-based/random approaches tend to mispredict the best first experiment, resulting in a substantial performance gap when compared to the MOCU-based approaches. Figure~\ref{fig10} also shows that the predicted experimental sequences diverge in later iterations. However, this does not impact the OED performance on average, as later experiments do not reduce the model uncertainty as significantly as the earlier experiments.

\subsubsection{Optimal experimental design for uncertain Kuramoto models with seven oscillators}
We also repeated the experiments for uncertain Kuramoto models that consist of seven oscillators. As before, true (unknown) models were randomly sampled from the uncertainty class one hundred times to evaluate average performance. We used the same parameters and model described in Sec.~\ref{MOCU7}. 

\begin{figure}
\centering
\includegraphics{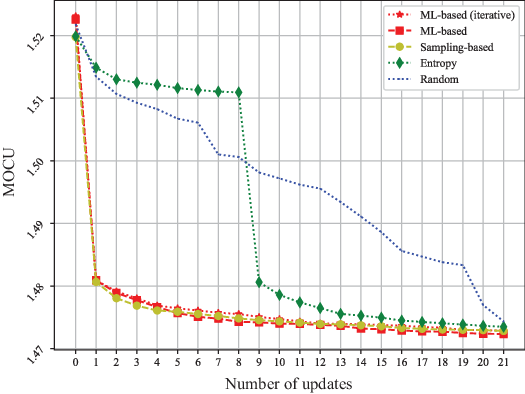}
\caption{Average performance of various experimental design strategies for uncertain Kuramoto models with seven oscillators based on one hundred experiments. All MOCU-based methods lead to the best performance, and random selection results in linear uncertainty reduction.}
\label{fig11}
\end{figure}

\begin{figure}
\centering
\includegraphics{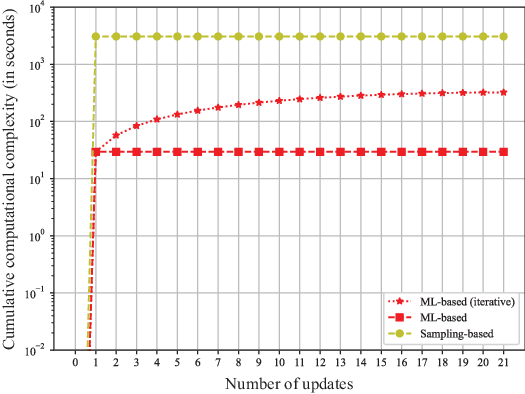}
\caption{
Average cumulative computational cost (in seconds) for the different experimental design strategies for uncertain Kuramoto models with seven oscillators.}
\label{fig12}
\end{figure}

Figure~\ref{fig11} shows the OED performance assessment results for the various experimental design methods based on the seven-oscillator Kuramoto model. As we can see from Fig.~\ref{fig11}, the performance trends were very similar to those seen in Fig.~\ref{fig8} for the Kuramoto model with five oscillators. The proposed ML-based methods again accurately identified the first optimal experiment that maximally reduces MOCU on average. All four MOCU-based OED schemes (both ML-based and sampling-based), regardless of whether or not the remaining expected MOCU values were re-estimated after each experimental update, showed almost identical performance on average. Figure~S3 compares the performance of different methods, where we measured the MOCU that remains after performing the first experiment selected by each method. The results are again shown for one hundred evaluations based on different true models. As shown in Fig.~S3, the efficacy of the first experiment varies depending on the underlying true model, which is expected. As before, the results in Fig.~S3 clearly show that the proposed ML-based OED scheme can effectively replicate the performance of the original sampling-based approach~\cite{Hong2021}, the primary goal of this study. The computational time is shown in Fig.~\ref{fig12}, which clearly shows that the ML-based scheme (especially, the non-iterative scheme) is significantly faster compared to the original sampling-based approach. Furthermore, even the ML-based method with the iterative update was considerably faster than the sampling-based that does not iteratively re-estimate the expected remaining MOCU. 

\begin{figure}
\centering
\includegraphics{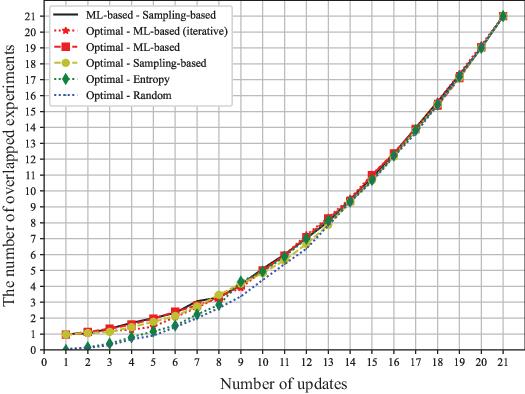}
\caption{
Comparison between the optimal sequence of experiments predicted by different OED strategies for uncertain Kuramoto model with seven oscillators. The $y$-axis shows the number of common experiments within the first $k$ experiments predicted by two different methods. The results are shown for uncertain Kuramoto models with seven oscillators.}
\label{fig13}
\end{figure}

Finally, we compared the experimental sequences identified by the different methods, including the true optimal experimental selection (\textit{i.e.}, predicted by an ``oracle''). Note that due to the excessive computational cost of the optimal experimental selection (as it requires exhaustive search), we identified the optimal experimental sequences only for the first thirty evaluations based on randomly sampled true models from the uncertainty class. For this reason, Fig.~\ref{fig13} shows the comparison results based on the first thirty experimental sequences (out of one hundred). As shown in the figure, both the ML-based and the sampling-based methods were able to accurately identify the first optimal experiment. The predicted sequences tend to diverge in later iterations. However, considering the simulation results shown in Fig.~\ref{fig11}, it is likely that this is because many experiments in later updates do not significantly reduce the objective uncertainty, once the best experiment has been performed in the earlier iterations (especially, the first iteration). Also, we can see that the entropy-based and the random selection approaches tend to miss the best experiment, which results in a significant degradation in the overall experimental design performance. These comprehensive simulation results clearly show that our proposed ML-based OED approach effectively quantifies the objective model uncertainty at a small fraction of the computational cost of the sampling-based method, thereby remarkably accelerating the OED process while maintaining excellent performance.

\section{Concluding Remarks}
\label{sec:conclusion}
In this paper, we proposed an ML approach that can significantly accelerate the objective-based quantification of model uncertainty via MOCU. A major bottleneck in applying MOCU for designing/prioritizing optimal experiments that can optimally reduce the uncertainty in models that represent real-world complex uncertain systems has been the high computational cost for accurately estimating MOCU. The proposed approach effectively addresses this issue in the context of OED for uncertain Kuramoto models by replacing the computational costly DE solver with an ML model, which remarkably speeds up the process of predicting the optimal controller (\textit{i.e.}, the oscillator that guarantees global frequency synchronization at minimum cost). The trained ML model predicts the asymptotic behavior of a given Kuramoto model, namely, whether all oscillators in the model will be eventually frequency synchronized or not. 

Note that there may be an outward similarity between the proposed ML-based approach in the context of the uncertain Kuramoto model with other earlier studies~\cite{L712178, Han8505, RAISSI2019686} whose main interest is to efficiently solve DEs via ML models,  as the prediction result of our trained ML model includes the process of solving Kuramoto model equations, albeit implicitly. However, more precisely, the trained ML model in the current study makes a ``decision'' based on the scientific knowledge regarding the underlying Kuramoto model for the purpose of identifying the optimal experiment that is expected to be most effective in reducing model uncertainty, which pertains to a broader issue rather than solving DEs as in other previous studies.

The results shown in Sec.~\ref{sec:mocu_ml} and Sec.~\ref{sec:oed_ml} clearly demonstrate that the ML-based MOCU calculations are highly correlated with those computed by the sampling-based scheme originally proposed in~\cite{Hong2021}. Furthermore, the OED performance of the ML-based scheme is practically equivalent to that of the original sampling-based OED scheme. However, despite achieving equivalent OED performance, our proposed ML-based OED scheme accelerates the experimental design process by at least two orders of magnitude, resulting in significant computational gains. The remarkably enhanced computational efficiency enables more reliable MOCU calculation by further increasing the sample size (\textit{i.e.}, $K$) as needed. Furthermore, it allows us to iteratively recompute the remaining MOCU $R\left(i,j\right)$ after performing the predicted optimal experiment at each experimental update (see Fig.~\ref{fig9} and Fig.~\ref{fig12}), which can--in theory--lead to a more accurate prediction of the optimal experiment, although the actual gain will depend on the underlying model uncertainty. Such iterative update is practically infeasible for the original sampling-based OED scheme without resorting to HPC (high-performance computing).

As shown in Sec.~\ref{sec:mocu_ml} and Sec.~\ref{sec:oed_ml}, our ML-based MOCU estimation and OED approach remarkably enhance the computational efficiency by refraining from repeatedly solving the DEs for the uncertain Kuramoto models for the sake of finding the optimal robust operator (which is required in the original sampling-based approach) but instead adopting ML for decision making. However, as training the ML model requires the generation of sufficient training data, which also requires solving the coupled ODEs for different Kuramoto models in the uncertainty class, it will be interesting to compare the proposed ML-based approach with the sampling-based approach from the perspective of ``data efficiency''. For this purpose, we quantitatively compare the proposed approach with the sampling-based approach in terms of data requirements. For the uncertain Kuramoto model with five oscillators, we trained the ML model (an fcNN with a single hidden layer) with $40,000$ labeled sample points. Each sample point corresponds to the Kuramoto model with a different parameter, and labeling the sample point (\textit{i.e.}, synchronized vs non-synchronized) requires solving the corresponding ODEs. The trained model is used throughout the entire experimental design process without the need for generating additional sample points. On the other hand, the sampling-based method requires generating approximately $2.2\times 10^7$ labeled sample points (\textit{i.e.}, by solving the DEs for different Kuramoto model parameters). Similarly, for the uncertain Kuramoto model with seven oscillators, we trained an fcNN model based on $100,000$ labeled sample points, and the trained model is used throughout the experimental design process. In comparison, the sampling-based approach requires the generation of around $9.4\times 10^7$ labeled sample points. These comparisons clearly show that our proposed ML-based OED acceleration scheme not only improves the computational efficiency but also drastically improves the data efficiency.

It is worth noting that the proposed approach, applying ML models to the estimation of the (remaining expected) MOCU can be generalized, extended, and applied to other MOCU-based OED problems concerning real-world applications that do not possess closed-form (remaining expected) MOCU. In such cases, based on the applications and the data types, one may consider different ML models including convolutional neural networks (CNN)~\cite{fukushima1982neocognitron}, recurrent neural networks (RNN)~\cite{rumelhart1986learning}, long short term memory networks (LSTM)~\cite{hochreiter1997long}, and graph convolutional networks (GCN)~\cite{scarselli2008graph}. While we have focused on accelerating the experimental design for uncertain Kuramoto models using ML, the proposed ML-based OED acceleration scheme is fairly general and its applicability goes beyond the Kuramoto models. For example, the proposed scheme may be used to accelerate the design of effective experiments to reduce uncertainty and improve the control performance of various other engineering models in the presence of uncertainty~\cite{tao2020robust,tao2021robust,zhang2021asynchronous,xin2022online}.

An interesting direction for future research is to utilize ML models to learn scientific knowledge from data. In this paper, we considered pairwise synchronization experiments, whose result can be used to reduce model uncertainty using Theorem~\ref{theo:sync} that gives us the necessary and sufficient condition for frequency synchronization of an oscillator pair. In the absence of such knowledge, one cannot design experiments for effective uncertainty reduction. As mentioned before, similar theorems do not exist in general for non-homogeneous Kuramoto models that consist of more than two oscillators. Discovering useful relational knowledge regarding the model parameters via ML can lead to the design of more effective experiments as well as a significant expansion of the potential experimental design space. We are currently investigating the potential utilization of deep neural network (DNN) models for knowledge discovery in non-homogeneous Kuramoto models with multiple oscillators.

\bibliographystyle{unsrt}  


\end{document}